\documentclass[10pt]{article}

\newcommand{\vekk}[1]{}

\usepackage[a4paper,text={140mm,222mm}]{geometry}
\usepackage[T1]{fontenc}
\usepackage[utf8]{inputenc}
\usepackage{amsmath,amssymb,amsthm,mathtools}
\usepackage{titlesec}

\usepackage[backend=biber, style=authoryear, natbib=true,
            giveninits=true, maxbibnames=99,
            doi=false, url=false, eprint=false]{biblatex}
\newbibmacro{string+doiurl}[1]{%
  \iffieldundef{doi}
    {\iffieldundef{url}
       {#1}
       {\href{\thefield{url}}{#1}}}
    {\href{https://doi.org/\thefield{doi}}{#1}}}
\DeclareFieldFormat{title}{\usebibmacro{string+doiurl}{\mkbibemph{#1}}}
\DeclareFieldFormat[article,inbook,incollection,inproceedings,patent,thesis,unpublished]{title}%
  {\usebibmacro{string+doiurl}{\mkbibquote{#1\isdot}}}

\usepackage[hidelinks]{hyperref}

\usepackage{orcidlink}
\usepackage{doi}

\titleformat{\section}{\centering\normalfont\scshape}{\thesection.}{0.5em}{}
\titlespacing*{\section}{0pt}{1.6ex plus .4ex minus .2ex}{0.9ex}
\vekk{
\pagestyle{fancy}\fancyhf{}
\fancyhead[LE]{\small\thepage\hfill\scshape Gunnar Taraldsen}
\fancyhead[RO]{\small\scshape Behrens-Fisher problem\hfill\thepage}

}

\newtheoremstyle{bika}{0.9ex}{0.9ex}{\itshape}{}{\scshape}{.}{ }{}
\theoremstyle{bika}
\newtheorem{theorem}{Theorem}
\newtheorem{lemma}{Lemma}
\newtheorem{proposition}{Proposition}
\theoremstyle{remark}
\newtheorem{example}{Example}
\renewenvironment{proof}[1][Proof]{\par\noindent\textit{#1}.\ }{\hfill$\square$\par\smallskip}

\DeclarePairedDelimiter{\abs}{\lvert}{\rvert}
\DeclareMathOperator{\pr}{pr}
\newcommand{\R}{\mathbb{R}}
\newcommand{\E}{\mathbb{E}}
\newcommand{\sgn}{\operatorname{sgn}}

\begin{document}

\begin{center}
{\large Validity of the Behrens--Fisher solution}\\[1.2ex]
\doi{10.48550/arXiv.2611.06987} \\[1.2ex]
{\scshape By Gunnar Taraldsen\orcidlink{0000-0003-4980-7019}} 
\vekk{{\itshape Department of Mathematical Sciences, Norwegian University of Science and
Technology, 7491 Trondheim, Norway}\\
\texttt{gunnar.taraldsen@ntnu.no}\\[1.6ex]  
}
\end{center}
\noindent
The Behrens--Fisher fiducial interval for the difference of two normal means is known
from the numerical work of \citet{Robinson76} to be conservative, but a proof has been
lacking. The results here prove that the fiducial interval for any linear combination of $k$ normal
means has coverage probability at least its nominal level for all values of the
unknown variances. 
\citet{Taraldsen26BehrensFisherproblemsolved} explains this in more detail
using the original data of \citet{Behrens29}.
The proof presented in the following rests on a novel plug-in inequality: if a threshold function
controls a scaled sum of pivots pointwise in the scale vector, then it also controls
the sum with the observed scales plugged in, provided the location errors satisfy a
concentration condition formulated through a weighted geometric mean. Independent
normal errors satisfy the condition.

\medskip\noindent
{\itshape Some key words}: Behrens--Fisher problem; Concentration inequality;
Confidence interval; Fiducial inference; Measurement uncertainty; Peakedness.

\section{Introduction}

\citet{Behrens29} compared the yields of grain varieties from small samples with
unequal variances, and \citet{Fisher35} gave the fiducial solution now known as the
Behrens--Fisher solution. Whether the resulting interval for the difference of two
normal means has frequentist coverage at least its nominal level, uniformly in the
unknown variance ratio, was investigated numerically by \citet{Robinson76}, who
reported that it does; a mathematical proof has not been available
\citep[p.~517]{LehmannRomano22}. This note gives a proof, and more generally proves the
validity of the fiducial interval for a linear combination of $k$ normal means
\citep{James59}, the case that arises in the evaluation of measurement uncertainty
\citep{JCGM08}.
The argument rests on a plug-in inequality for scaled sums, Theorem~\ref{thm:main},
whose only hypothesis beyond independence is a concentration condition on the
location errors formulated through a weighted geometric mean. Independent normal
errors satisfy the condition, Proposition~\ref{prop:normal}.

\section{The plug-in inequality and the concentration condition}
\label{sec:ineq}

Let $X=(X_1,\ldots,X_k)$ and $S=(S_1,\ldots,S_k)$ be independent random vectors with
$S_i>0$. The $X_i$ take values in a separable Banach space $E$ with norm
$\abs{\cdot}$; the reader may take $E=\R$ throughout, since only the measurability and
the positive homogeneity $\abs{ax}=a\abs{x}$, $a>0$, of the norm are used. Fix
$\varepsilon\in(0,1)$ and let $d$ be a measurable function on $(0,\infty)^k$ such
that
\begin{equation}\label{eq:pointwise}
\pr\Big(\Big\lvert\sum_{i=1}^k s_i\frac{X_i}{S_i}\Big\rvert\le d(s)\Big)\ge 1-\varepsilon
\qquad\text{for all } s\in(0,\infty)^k .
\end{equation}
Since \eqref{eq:pointwise} holds for every $s$, it is tempting 
to substitute $S$ for $s$ and conclude that
\begin{equation}\label{eq:plugin}
\pr\Big(\Big\lvert\sum_{i=1}^k X_i\Big\rvert\le d(S)\Big)\ge 1-\varepsilon .
\end{equation}
This 'proof' of \eqref{eq:plugin} is wrong.
In the Behrens--Fisher problem the $X_i/S_i$ are Student pivots,
\eqref{eq:pointwise} with the smallest admissible $d$ defines the fiducial interval,
and \eqref{eq:plugin} is its validity as explained by \citet{Taraldsen26BehrensFisherproblemsolved}.
The smallest admissible
threshold is
\begin{equation}\label{eq:dstar}
d^\ast(s)=\inf\Big\{t\ge 0:\ \pr\Big(\Big\lvert\sum_{i=1}^k s_i\frac{X_i}{S_i}\Big\rvert\le t\Big)\ge 1-\varepsilon\Big\},
\end{equation}
whose properties are given in Lemma~\ref{lem:dstar} below. 

\vekk{
\begin{example}\label{ex:counter}
Let $k=2$, $X_1$ uniform on $\{-1,0,1\}$, $X_2$ uniform on $\{-2,0,2\}$, $S$ uniform
on $\{(2,1),(1,1)\}$, all independent, $E=\R$ and $\varepsilon=4/5$. Given
$S=(1,1)$ the sum in \eqref{eq:pointwise} at $s=(2,1)$ is $2X_1+X_2$, which vanishes
with probability $1/3$, so $\pr(\lvert\sum_i s_iX_i/S_i\rvert\le 0)=
\tfrac12\cdot\tfrac19+\tfrac12\cdot\tfrac13=2/9\ge 1/5$ and $d^\ast(2,1)=0$.
Similarly $d^\ast(1,1)=1/2$, since $\lvert X_1/2+X_2\rvert\le 1/2$ has probability
$1/3$ and $\pr(\lvert X_1+X_2\rvert\le 1/2)=1/9$. Thus $d^\ast$ satisfies
\eqref{eq:pointwise}, but the left side of \eqref{eq:plugin} is
$\tfrac12\pr(\lvert X_1+X_2\rvert\le 0)+\tfrac12\pr(\lvert X_1+X_2\rvert\le\tfrac12)
=1/9<1/5$.
\end{example}
}

A random vector $V$ in $E$ is \emph{more concentrated} about $0$ than $W$ if
\begin{equation}\label{eq:peak}
\pr(\abs{V}\le t)\ge\pr(\abs{W}\le t)\qquad\text{for all } t\ge 0 ,
\end{equation}
which is the peakedness ordering of \citet{Birnbaum48}. The vector $X$ satisfies the
\emph{concentration condition} if there is a probability vector $w=(w_1,\ldots,w_k)$
such that $\tilde c\sum_iX_i$ is more concentrated than $\sum_ic_iX_i$ for all
$c\in(0,\infty)^k$, where $\tilde c=\prod_ic_i^{w_i}$ is the weighted geometric mean.
Replacing $\tilde c$ by the weighted arithmetic mean $\bar c=\sum_iw_ic_i$ gives a
stronger condition, since $\tilde c\le\bar c$ and $a\mapsto\pr(\abs{a\sum_iX_i}\le t)$
is nonincreasing. 
\vekk{The vector in Example~\ref{ex:counter} violates the condition for
every $w$: with $c=(2,1)$ one has $\tilde c=2^{w_1}\ge 1$, so
$\pr(\lvert\tilde c(X_1+X_2)\rvert\le t)=1/9<1/3=\pr(\lvert 2X_1+X_2\rvert\le t)$ for
$0\le t<1$.}

\begin{proposition}\label{prop:normal}
Let $X_i=\sigma_iZ_i$, where $\sigma_i>0$ and $Z_1,\ldots,Z_k$ are independent
standard normal variables. Then $X$ satisfies the concentration condition with
$w_i=\sigma_i^2/\sum_j\sigma_j^2$.
\end{proposition}
\begin{proof}
Both $\sum_i c_i X_i$ and $\bar c \sum_iX_i$ are centred normal, with variances
$\sum_i c_i^2\sigma_i^2 \ge \bar c^{\,2} \sum_j\sigma_j^2$ from
Jensen's inequality.
The conclusion folows since a 
centred normal variable is more concentrated than another if its variance is smaller.
\end{proof}
The proof applies verbatim to $X_i=\sigma_iZ_i$ with $Z_1,\ldots,Z_k$ independent
copies of a centred Gaussian vector $Z$ in $E$, since $\sum_ia_iZ_i$ then has the law
of $(\sum_ia_i^2)^{1/2}Z$, and to independent symmetric $\alpha$-stable $Z_i$, $0<\alpha\le 2$, with
$w_i\propto\sigma_i^\alpha$, by the inequality between the weighted power means of
orders $\alpha$ and $0$.

\section{The main result}
\label{sec:main}

\begin{lemma}\label{lem:dstar}
Let $T=(T_1,\ldots,T_k)$ be a random vector in $E^k$ and $a\in(0,1)$. The function
\begin{equation}\label{eq:quantile}
d^\ast(s)=\inf\Big\{t\ge 0:\ \pr\Big(\Big\lvert\sum_{i=1}^k s_iT_i\Big\rvert\le t\Big)\ge a\Big\},
\qquad s\in(0,\infty)^k,
\end{equation}
is finite and measurable, the infimum is attained, and $d^\ast$ is positively
homogeneous: $d^\ast(cs)=c\,d^\ast(s)$ for $c>0$.
\end{lemma}

\begin{proof}
Put $F_s(t)=\pr(\lvert\sum_is_iT_i\rvert\le t)$. For fixed $s$ this is
nondecreasing in $t$, tends to $1$ as $t\to\infty$, and is right-continuous, since
the events $\{\lvert\sum_is_iT_i\rvert\le t\}$ decrease to
$\{\lvert\sum_is_iT_i\rvert\le t_0\}$ as $t\downarrow t_0$. Hence the set in
\eqref{eq:quantile} is a nonempty interval $[d^\ast(s),\infty)$, and the infimum is
attained. Homogeneity follows by the substitution $t=cu$ in \eqref{eq:quantile},
because $\lvert\sum_ics_iT_i\rvert=c\lvert\sum_is_iT_i\rvert$. 

For measurability, note that the indicator function
$(s,t,\omega)\mapsto (\lvert\sum_is_iT_i(\omega)\rvert \le t)$ 
is measurable, so that
$(s,t)\mapsto F_s (t)$ is measurable by Fubini's theorem
\citep[Theorem~8.8]{Rudin87} and $A=\{(s,t):F_s(t)\ge a\}$ is a measurable set.
Since the infimum is attained, $\{s:d^\ast(s)\le t\}=\{s:F_s(t)\ge a\}$ is the
$t$-section of $A$, which is measurable \citep[Theorem~8.2]{Rudin87}.
\end{proof}

\begin{theorem}\label{thm:main}
Let $X$ satisfy the concentration condition with probability weight $w$, let $S>0$ be
independent of $X$, and let $\varepsilon\in(0,1)$. Then \eqref{eq:pointwise} implies
\eqref{eq:plugin} for every measurable $d$.
\end{theorem}

\begin{proof}
Write $\tilde s=\prod_is_i^{w_i}$ and $\tilde S=\prod_iS_i^{w_i}$, and let $v$ denote
a vector with $\tilde v=1$.

\textit{Step 1: homogeneous $d$.} Assume $d(as)=a\,d(s)$ for $a>0$,
so \eqref{eq:pointwise} gives
\begin{equation}\label{eq:step1}
\pr\Big(\Big\lvert\sum_i\frac{s_i}{\tilde s}\,\frac{\tilde S}{S_i}X_i\Big\rvert
\le\tilde S\,d(\frac{s}{\tilde s})\Big)\ge 1-\varepsilon .
\end{equation}
Condition on $S=u$. The coefficient vector $c_i=(s_i/\tilde s)(\tilde u/u_i)$ has
$\tilde c=1$, so the concentration condition and the independence of $X$ and $S$ give
$\pr(\lvert\sum_ic_iX_i\rvert\le\tilde u\,d(s/\tilde s))\le
\pr(\lvert\sum_iX_i\rvert\le\tilde u\,d(s/\tilde s))$. Integrating over the law of
$S$ and using \eqref{eq:step1} yields, as $v=s/\tilde s$ ranges over all vectors
with $\tilde v=1$,
\begin{equation}\label{eq:step2}
\pr\Big(\Big\lvert\sum_iX_i\Big\rvert\le\tilde S\,d(v)\Big)\ge 1-\varepsilon
\qquad\text{for all } v .
\end{equation}
Now let $\delta_\ast=\inf_vd(v)\ge 0$ and choose $v_n$ with $d(v_n)\downarrow\delta_\ast$.
Homogeneity gives $d(S)=\tilde S\,d(S/\tilde S)\ge\tilde S\delta_\ast$, and the
events $\{\lvert\sum_iX_i\rvert\le\tilde Sd(v_n)\}$ decrease to
$\{\lvert\sum_iX_i\rvert\le\tilde S\delta_\ast\}$. Continuity of probability and
\eqref{eq:step2} therefore give
\begin{equation}\label{eq:step3}
\pr\Big(\Big\lvert\sum_iX_i\Big\rvert\le d(S)\Big)
\ge\pr\Big(\Big\lvert\sum_iX_i\Big\rvert\le\tilde S\delta_\ast\Big)
=\lim_{n}\pr\Big(\Big\lvert\sum_iX_i\Big\rvert\le\tilde S\,d(v_n)\Big)\ge 1-\varepsilon .
\end{equation}

\textit{Step 2: general $d$.} Apply Lemma~\ref{lem:dstar} with $T_i=X_i/S_i$ and
$a=1-\varepsilon$, so that $d^\ast$ is the function \eqref{eq:dstar}. Since the
infimum in \eqref{eq:quantile} is attained, $d^\ast$ satisfies
\eqref{eq:pointwise}, and every measurable $d$ satisfying \eqref{eq:pointwise}
obeys $d\ge d^\ast$. By Lemma~\ref{lem:dstar}, $d^\ast$ is
measurable and homogeneous, so Step 1 gives \eqref{eq:plugin} for $d^\ast$, and
\eqref{eq:plugin} for $d$ follows from
$\{\lvert\sum_iX_i\rvert\le d^\ast(S)\}\subseteq\{\lvert\sum_iX_i\rvert\le d(S)\}$.
\end{proof}

The proof uses the norm only through its positive homogeneity and measurability,
so Theorem~\ref{thm:main} holds for any measurable positively homogeneous
$\abs{\cdot}\ge 0$, symmetric or not. The intermediate bound in \eqref{eq:step3}
shows that in the normal model the inequality in \eqref{eq:plugin} is strict whenever
$d(S/\tilde S)>\delta_\ast$ with positive probability, which accounts for the strict
conservativeness observed numerically by \citet{Robinson76}.

\printbibliography[title={References}]

\end{document}